# Optimum Statistical Test Procedure

Rajesh Singh

Department of Statistics, BHU, Varanasi (U.P.), India

rsinghstat@yahoo.com

Jayant Singh

Department of Statistics

Rajasthan University, Jaipur, India

Jayantsingh47@rediffmail.com

Florentin Smarandache
Chair of Department of Mathematics, University of New Mexico, Gallup, USA

smarand@unm.edu

## Introduction

Let X be a random variable having probability distribution $P(X/\theta)$, $\theta \in \Theta$. It is desired to test $H_0 : \theta \in \Theta_0$ against $H_1 : \theta \in \Theta_1 = \Theta - \Theta_0$. Let S denote the sample space of outcomes of an experiment and $\underline{x} = (x_1, x_2, ---, x_n)$ denote an arbitrary element of S. A test procedure consists in diving the sample space into two regions W and $S - W$ and deciding to reject $H_0$ if the observed $\underline{x} \in W$. The region W is called the critical region. The function $\gamma(\theta) = P_\theta(\underline{x} \in W) = P_\theta(W)$, say, is called the power function of the test.

We consider first the case where $\Theta_0$ consists of a single element, $\theta_0$ and its complement $\Theta_1$ also has a single element $\theta_1$. We want to test the simple hypothesis $H_0 : \theta = \theta_0$ against the simple alternative hypothesis $H_1 : \theta = \theta_1$.

Let $L_0 = L(X/H_0)$ and $L_1 = L(X/H_1)$ be the likelihood functions under $H_0$ and $H_1$ respectively. In the Neyman – Pearson set up the problem is to determine a critical region W such that

$$\gamma(\theta_0) = P_{\theta_0}(W) = \int_W L_0 dx = \alpha, \text{ an assigned value} \tag{1}$$

and $$\gamma(\theta_1) = P_{\theta_1}(W) = \int_W L_1 dx \quad \text{is maximum} \tag{2}$$

compared to all other critical regions satisfying (1).

If such a critical region exists it is called the most powerful critical region of level $\alpha$.

By the Neyman-Pearson lemma the most powerful critical region $W_0$ for testing $H_0\colon \theta = \theta_0$ against $H_1\colon \theta = \theta_1$ is given by

$$W_0 = \{\underline{x} : L_1 \geq kL_0\}$$

provided there exists a k such that (1) is satisfied.

For this test $\gamma(\theta_0) = \alpha$ and $\gamma(\theta_1) \to 1$ as $n \to \alpha$.

But for any good test we must have $\gamma(\theta_0) \to 0$ and $\gamma(\theta_1) \to 1$ as $n \to \infty$ because complete discrimination between the hypotheses $H_0$ and $H_1$ should be possible as the sample size becomes indefinitely large.

Thus for a good test it is required that the two error probabilities $\alpha$ and $\beta$ should depend on the sample size n and both should tend to zero as $n \to \infty$.

We describe below test procedures which are optimum in the sense that they minimize the sum of the two error probabilities as compared to any other test. Note that minimizing $(\alpha+\beta)$ is equivalent to maximising

$$1 - (\alpha+\beta) = (1-\beta)-\alpha = \text{Power} - \text{Size}.$$

Thus an optimum test maximises the difference of power and size as compared to any other test.

**Definition 1** : A critical region $W_0$ will be called optimum if

$$\int_{W_0} L_1 dx - \int_{W_0} L_0 dx \geq \int_{W} L_1 dx - \int_{W} L_0 dx \qquad (3)$$

for every other critical region W.

**Lemma 1**: For testing $H_0$: $\theta = \theta_0$ against $H_1$ : $\theta = \theta_1$ the region

$W_0 = \{\underline{x} : L_1 \geq L_0\}$ is optimum.

**Proof** : $W_0$ is such that inside $W_0$, $L_1 \geq L_0$ and outside $W_0$ , $L_1 < L_0$. Let W be any other critical region.

$$\int_{W_0} (L_1 - L_0) dx - \int_{W} (L_1 - L_0) dx$$

$$= \int_{W_0 \cap W^c} (L_1 - L_0) dx - \int_{W \cap W_0^c} (L_1 - L_0) dx,$$

by subtracting the integrals over the common region $W_0 \cap W$.

$\geq 0$

since the integrand of first integral is positive and the integrand under second integral is negative.

Hence (3) is satisfied and $W_0$ is an optimum critical region.

**Example 1** :

Consider a normal population $N(\theta, \sigma^2)$ where $\sigma^2$ is known.

It is desired to test $H_0$: $\theta = \theta_0$ against $H_1$ : $\theta = \theta_1$ , $\theta_1 > \theta_0$ .

$$L(x/\theta) = \left(\frac{1}{\sigma\sqrt{2\pi}}\right)^n e^{-\sum_{i=1}^{n}\frac{(x_i-\theta)^2}{2\sigma^2}}$$

$$\frac{L_1}{L_0} = \frac{e^{-\frac{\sum(x_i-\theta_1)^2}{2\sigma^2}}}{e^{-\frac{\sum(x_i-\theta_0)^2}{2\sigma^2}}}$$

The optimum test rejects $H_o$

if $\quad \frac{L_1}{L_0} \geq 1$

i.e. if $\quad \log\frac{L_1}{L_0} \geq 0$

i.e. if $\quad -\frac{\sum(x_i-\theta_1)^2}{2\sigma^2} + \frac{\sum(x_i-\theta_0)^2}{2\sigma^2} \geq 0$

i.e. if $\quad 2\theta_1\sum x_i - n\theta_1^2 - 2\mu_0\sum x_i + n\theta_0^2 \geq 0$

i.e. if $\quad (2\theta_1 - 2\theta_0)\sum x_i \geq n\left(\theta_1^2 - \theta_0^2\right)$

i.e. if $\dfrac{\sum x_i}{n} \geq \dfrac{\theta_1+\theta_0}{2}$

i.e. if $\bar{x} \geq \dfrac{\theta_1+\theta_0}{2}$

Thus the optimum test rejects $H_0$

if $\bar{x} \geq \dfrac{\theta_1+\theta_0}{2}$

We have

$$\alpha = P_{H_0}\left[\bar{x} \geq \frac{\theta_1+\theta_0}{2}\right]$$

$$= P_{H_0}\left[\frac{\bar{x}-\theta_0}{\sigma/\sqrt{n}} \geq \frac{\sqrt{n}(\theta_1-\theta_0)}{2\sigma}\right]$$

Under $H_o$, $\dfrac{\bar{x}-\theta_0}{\left(\sigma/\sqrt{n}\right)}$ follows N(0,1) distribution.

$$\therefore \quad \alpha = 1-\Phi\left(\frac{\sqrt{n}(\theta_1-\theta_0)}{2\sigma}\right)$$

where $\Phi$ is the c.d.f. of a N(0,1) distribution.

$$\beta = P_{H_1}\left[\bar{x} < \frac{\theta_1+\theta_0}{2}\right] = P_{H_1}\left[\frac{\bar{x}-\theta_1}{\sigma/\sqrt{n}} < \frac{-\sqrt{n}(\theta_1-\theta_0)}{2\sigma}\right]$$

Under $H_1$, $\dfrac{\bar{x}-\theta_1}{\left(\sigma/\sqrt{n}\right)}$ follows N(0,1) distribution.

$$\beta = 1 - \Phi\left(\frac{\sqrt{n}(\theta_1 - \theta_0)}{2\sigma}\right)$$

Here $\alpha = \beta$.

It can be seen that $\alpha = \beta \to 0$ as $n \to \infty$.

**Example 2** : For testing $H_0$: $\theta = \theta_0$ against $H_1$ : $\theta = \theta_1$ , $\theta_1 < \theta_0$, the optimum test rejects

$H_0$ when $\bar{x} \le \frac{\theta_1 + \theta_0}{2}$.

**Example 3** : Consider a normal distribution $N(\theta, \sigma^2)$, $\theta$ known.

It is desired to test $H_0$ : $\sigma^2 = \sigma_0^2$ against H1 : $\sigma^2 = \sigma_1^2, \sigma_1^2 > \sigma_0^2$.

We have

$$L(x/\sigma^2) = \left(\frac{1}{2\pi\sigma^2}\right)^{\frac{n}{2}} e^{-\sum \frac{(x_i - \theta)^2}{2\sigma^2}}$$

$$\frac{L_1}{L_0} = \frac{\left(\frac{\sigma_0^2}{\sigma_1^2}\right)^{\frac{n}{2}} e^{-\frac{\sum (x_i - \theta)^2}{2\sigma_1^2}}}{e^{-\frac{\sum (x_i - \theta)^2}{2\sigma_0^2}}}$$

$$\log \frac{L_1}{L_0} = -\frac{n}{2}\left(\log \sigma_1^2 - \log \sigma_0^2\right) - \frac{\sum (x_i - \theta)^2}{2\sigma_1^2} + \frac{\sum (x_i - \theta)^2}{2\sigma_0^2}$$

$$= -\frac{n}{2}\left(\log \sigma_1^2 - \log \sigma_0^2\right) + \frac{\sum (x_i - \theta)^2}{2} \frac{\sigma_1^2 - \sigma_0^2}{\sigma_1^2 \sigma_0^2}$$

The optimum test rejects $H_0$

$$if \quad \frac{L_1}{L_0} \geq 1$$

i.e. if $$\frac{\sigma_1^2 - \sigma_0^2}{2\sigma_1^2\sigma_0^2}\sum(x_i - \theta)^2 \geq \frac{n}{2}\left(\log\sigma_1^2 - \log\sigma_0^2\right)$$

i.e. if $$\frac{\sum(x_i - \theta)^2}{\sigma_0^2} \geq \frac{n\sigma_1^2}{\sigma_1^2 - \sigma_0^2}\left(\log\sigma_1^2 - \log\sigma_0^2\right)$$

i.e. if $$\sum\left(\frac{x_i - \theta}{\sigma_0}\right)^2 \geq nc$$

where $$c = \frac{\sigma_1^2}{\sigma_1^2 - \sigma_0^2}\left(\log\sigma_1^2 - \log\sigma_0^2\right)$$

Thus the optimum test rejects $H_0$ if $\sum\left(\frac{x_i - \theta}{\sigma_0}\right)^2 \geq \ nc.$

Note that under $H_0 : \frac{x_i - \theta}{\sigma_0}$ follows N(0,1) distribution. Hence $\sum\left(\frac{x_i - \theta}{\sigma_0}\right)^2$ follows, under $H_0$, a chi-square distribution with n degrees of freedom (d.f.).

Here $$\alpha = P_{H_0}\left[\sum\left(\frac{x_i - \theta_0}{\sigma_0}\right)^2 \geq nc\right] = P\left[\chi^2_{(n)} \geq nc\right]$$

and $$1 - \beta = P_{H_1}\left[\sum\left(\frac{x_i - \theta}{\sigma_0}\right)^2 \geq nc\right]$$

$$= P_{H_1}\left[\sum\left(\frac{x_i - \theta}{\sigma_1}\right)^2 \geq \frac{nc\sigma_0^2}{\sigma_1^2}\right]$$

$$= P_{H_1}\left[\chi^2_{(n)} \geq \frac{nc\sigma_0^2}{\sigma_1^2}\right]$$

Note that under H1, $\sum\left(\frac{x_i-\theta}{\sigma_1}\right)^2$ follows a chi-square distribution with n d.f.

It can be seen that $\alpha \to 0$ and $\beta \to 0$ as $n \to \infty$.

**Example 4** : Let X follow the exponential family distributions

$$f(x/\theta) = c(\theta)e^{Q(\theta)T(x)}h(x)$$

It is desired to test $H_0$: $\theta = \theta_0$ against $H_1$ : $\theta = \theta_1$

$$L(x/\theta) = [c(\theta)]^n e^{Q(\theta)\sum T(x_i)} \prod_i h(x_i)$$

The optimum test rejects $H_o$ when

$$log\frac{L_1}{L_0} \geq 0$$

i.e. when $[Q(\theta_1) - Q(\theta_0)]\sum T(x_i) \geq n\log\frac{c(\theta_0)}{c(\theta_1)}$

i.e. when $\sum T(x_i) \geq \frac{n\log\frac{c(\theta_0)}{c(\theta_1)}}{[Q(\theta_1) - Q(\theta_0)]}$ if $Q(\theta_1) - Q(\theta_0) > 0$

and rejects $H_o$,

when $\sum T(x_i) \leq \frac{n\log\frac{c(\theta_0)}{c(\theta_1)}}{[Q(\theta_1) - Q(\theta_0)]}$ if $Q(\theta_1) - Q(\theta_0) < 0$

**Locally Optimum Tests:**

Let the random variable X have probability distribution $P(x/\theta)$. We are interested in testing $H_0: \theta = \theta_0$ against $H_1: \theta > \theta_0$. If W is any critical region then the power of the test as a function of $\theta$ is

$$\gamma(\theta) = P_\theta(W) = \int_W L(X/\theta)dx$$

We want to determine a region W for which

$$\gamma(\theta) - \gamma(\theta_0) = \int_W L(x/\theta)dx - \int_W L(x/\theta)$$

is a maximum.

When a uniformly optimum region does not exist, there is not a single region which is best for all alternatives. We may, however, find regions which are best for alternatives close to the null hypothesis and hope that such regions will also do well for distant alternatives. We shall call such regions locally optimum regions.

Let $\gamma(\theta)$ admit Taylor expansion about the point $\theta = \theta_0$. Then

$$\gamma(\theta) = \gamma(\theta_0) + (\theta - \theta_0)\acute{\gamma}(\theta_0) + \delta \qquad \text{where } \delta \to 0 \text{ as } \theta \to \theta_0$$

$$\therefore \quad \gamma(\theta) - \gamma(\theta_0) = (\theta - \theta_0)\acute{\gamma}(\theta_0) + \delta$$

If $|\theta - \theta_0|$ is small $\gamma(\theta) - \gamma(\theta_0)$ is maximised when $\acute{\gamma}(\theta_0)$is maximised.

**Definition2 :** A region $W_0$ will be called a locally optimum critical region if

$$\int_{W_0} \acute{L}(x/\theta_0)dx \geq \int_W \acute{L}(x/\theta_0)dx \tag{4}$$

For every other critical region W.

**Lemma 2** : Let $W_0$ be the region $\{\underline{x}: \acute{L}(x/\theta_0) \geq L(x/\theta_0)\}$. Then $W_0$ is locally optimum.

**Proof**: Let $W_0$ be the region such that inside it $\acute{L}(x/\theta_0) \geq L(x/\theta_0)$ and outside it

$\acute{L}(x/\theta_0) < L(x/\theta_0)$. Let W be any other region.

$$\int_{W_0} \acute{L}(x/\theta_0)dx - \int_{W} \acute{L}(x/\theta_0)dx$$

$$= \int_{W_{0\cap W^c}} \acute{L}(x/\theta_0)dx - \int_{W_0^c \cap W} \acute{L}(x/\theta_0)dx$$

$$= \int_{W_{0\cap W^c}} \acute{L}(x/\theta_0)dx + \int_{W_0^c \cup W} \acute{L}(x/\theta_0)dx \qquad (*)$$

$$\geq \int_{W_{0\cap W^c}} L(x/\theta_0)dx + \int_{W_0^c \cup W} L(x/\theta_0)dx$$

since $\acute{L}(x/\theta_0) \geq L(x/\theta_0)$ inside both the regions of the integrals.

$\geq 0$, since $L(x/\theta_0) \geq 0$ in all the regions.

Hence $\int_{W_0} \acute{L}(x/\theta_0)dx \geq \int_{W} \acute{L}(x/\theta_0)dx$ for every other region W.

**To prove (*):**

We have $\int_{R} L(x/\theta_0)dx + \int_{R^c} L(x/\theta_0)dx = 1$ for every region R.

Differentiating we have

$$\int_{R} \acute{L}(x/\theta_0)dx + \int_{R^c} \acute{L}(x/\theta_0)dx = 0$$

$$\int_{R} \acute{L}(x/\theta_0)dx = -\int_{R^c} \acute{L}(x/\theta_0)dx = 0$$

In (*), take $R^C = W_0^c \cap W$ and the relation is proved.

Similarly if the alternatives are $H_1 : \theta < \theta_0$ , the locally optimum critical region is

$\{\underline{x}: \acute{L}(x/\theta_0) \leq L(x/\theta_0)\}$.

**Example 5**: Consider $N(\theta, \sigma^2)$ distribution, $\sigma^2$ known.

It is desired to test $H_0$: $\theta = \theta_0$ against $H_1$ : $\theta > \theta_0$

$$L(x/\theta) = \left(\frac{1}{\sigma\sqrt{2\pi}}\right)^n e^{\frac{-\sum(x_i-\theta)^2}{2\sigma^2}}$$

$$logL(x/\theta) = nlog\left(\frac{1}{\sigma\sqrt{2\pi}}\right) - \frac{\sum(x_i - \theta)^2}{2\sigma^2}$$

$$\frac{\acute{L}(x/\theta)}{L(x/\theta)} = \frac{\delta logL(x/\theta)}{\delta\theta} = \frac{\sum(x_i - \theta)}{\sigma^2} = \frac{n(\bar{x} - \theta)}{\sigma^2}$$

$$\therefore \quad \frac{\acute{L}(x/\theta_0)}{L(x/\theta_0)} = \frac{n(\bar{x} - \theta_0)}{\sigma^2}$$

The locally optimum test rejects $H_0$ ,if

$$\frac{n(\bar{x} - \theta_0)}{\sigma^2} \geq 1$$

i.e. $\bar{x} \geq \theta_0 + \frac{\sigma^2}{n}$

Now,

$$\alpha \quad = P_{H_0}\left[\bar{x} \geq \theta_0 + \frac{\sigma^2}{n}\right]$$

$$= P_{H_0}\left[\frac{\bar{x} - \theta_0}{\sigma/\sqrt{n}} \geq \sigma/\sqrt{n}\right]$$

$= 1 - \Phi\left(\sigma/\sqrt{n}\right)$, since under $H_0$, $\frac{\bar{x}-\theta_0}{\sigma/\sqrt{n}}$ follows N(0,1) distribution.

$$1-\beta = P_{H_1}\left[\bar{x} \geq \theta_0 + \sigma^2/n\right]$$

$$= P_{H_1}\left[\frac{\bar{x}-\theta_1}{\sigma/\sqrt{n}} \geq -\frac{\theta_1-\theta_0}{\sigma/\sqrt{n}} + \frac{\sigma}{\sqrt{n}}\right]$$

$$= 1-\Phi\left[\frac{-(\theta_1-\theta_0\sqrt{n})}{\sigma} + \frac{\sigma}{\sqrt{n}}\right]$$

since under $H_1$, $\frac{\bar{x}-\theta_1}{\sigma/\sqrt{n}}$ follows N(0,1) distribution.

**Exercise** : If $\theta_0 = 10, \theta_1 = 11, \sigma = 2, n = 16$, then $\alpha = 0.3085$, $1\text{-}\beta = 0.9337$

Power - Size = 0.6252

If we reject $H_0$ when $\frac{\bar{x}-\theta_0}{\sigma/\sqrt{n}} > 1.64$, then $\alpha = 0.05$, $1\text{-}\beta = 0.6406$

Power – Size = 0.5906

Hence Power – Size of locally optimum test is greater than Power – size of the usual test.

**Locally Optimum Unbiased Test:**

Let the random variable X follows the probability distribution $P(x/\theta)$. Suppose it is desired to test $H_0$: $\theta = \theta_0$ against $H_1$ : $\theta \neq \theta_0$ . We impose the unbiasedness restriction

$\gamma(\theta) \geq \gamma(\theta_0), \theta \neq \theta_0$

and $\gamma(\theta) - \gamma(\theta_0)$ is a maximum as compared to all other regions. If such a region does not exist we impose the unbiasedness restriction $\acute{\gamma}(\theta_0) = 0$.

Let $\gamma(\theta)$ admit Taylor expansion about the point $\theta = \theta_0$.

Then $\gamma(\theta) = \gamma(\theta_0) + (\theta - \theta_0)\acute{\gamma}(\theta_0) + \frac{(\theta-\theta_0)^2}{2}\gamma''(\theta_0) + \eta$

where $\eta \to 0$ as $\theta \to 0$.

$\therefore \quad \gamma(\theta) - \gamma(\theta_0) = (\theta - \theta_0)\acute{\gamma}(\theta_0) + \frac{(\theta-\theta_0)^2}{2}\gamma''(\theta_0) + \eta$

Under the unbiasedness restriction $\acute{\gamma}(\theta_0) = 0$, if $|\theta - \theta_0|$ is small $\gamma(\theta) - \gamma(\theta_0)$ is maximised when $\gamma''(\theta_0)$ is maximised.

**Definition 3** : A region $W_0$ will be called a locally optimum unbiased region if

$$\acute{\gamma}(\theta_0) = \int_{W_0} \acute{L}(x/\theta_0)dx = 0 \tag{5}$$

and

$$\gamma''(\theta_0) = \int_{W_0} L''(X/\theta_0)dx \geq \int_{W} L''(X/\theta_0)dx \tag{6}$$

for all other regions W satisfying (5).

**Lemma 3** : Let $W_0$ be the region $\{\underline{x} : L''(x/\theta_0) \geq L(x/\theta_0)\}$

Then $W_0$ is locally optimum unbiased.

**Proof** : Let W be any other region

$$\int_{W_0} L''(x/\theta_0)dx - \int_{W} L''(x/\theta_0)dx$$

$$= \int_{W_0 \cap W^c} L''(x/\theta_0)dx + \int_{W_0^c \cap W} L''(x/\theta_0)dx$$

by subtracting the common area of W and $W_0$.

$$= \int_{W_0 \cap W^c} L''(x/\theta_0)dx + \int_{W_0 \cap W^c} L''(x/\theta_0)dx,$$

since $L''(x/\theta_0) \geq L(x/\theta)$ inside $W_0$ and outside W.

$\geq 0$ since $L(x/\theta) \geq 0.$

**Example 6**:

Consider $N(\theta, \sigma^2)$ distribution, $\sigma^2$ known .

$H_0: \theta = \theta_0, \ \ H_1: \theta \neq \theta_0$

$$L = \left(\frac{1}{\sigma\sqrt{2\pi}}\right)^n e^{-\frac{\sum(x_i-\theta)^2}{2\sigma^2}}$$

$$\frac{L"(x/\theta)}{L(x/\theta)} = \frac{n^2(\bar{x}-\theta)^2}{\sigma^4} - \frac{n}{\sigma^2}$$

Locally optimum unbiased test rejects $H_0$

$$if \ \ \frac{n^2(\bar{x}-\theta_0)^2}{\sigma^4} - \frac{n}{\sigma^2} \geq 1$$

$$i.e. \ \ \frac{n(\bar{x}-\theta_0)^2}{\sigma^2} \geq 1 + \frac{\sigma^2}{n}$$

$Under\ H_{0,}$ $\frac{n(\bar{x}-\theta_0)^2}{\sigma^2}$ follows $\chi^2_{(1)}$ distribution.

**Testing Mean of a normal population when variance is unknown.**

Consider $N(\theta, \sigma^2)$ distribution, $\sigma^2$ known.

For testing $H_0$: $\theta = \theta_0$ against $H_1$ : $\theta = \theta_1$, the critical function of the optimum test is given by

$$\phi_m(x) = \begin{cases} 1 & if\ L(x/\theta_1) \geq L(x/\theta_0) \\ 0 & otherwise \end{cases}$$

On simplification we get

$$\phi_m(x) = \begin{cases} 1 & if\ \ \bar{x} \geq \dfrac{\theta_0 + \theta_1}{2}\ \ if\ \mu_1 > \mu_0 \\ 0 & otherwise \end{cases}$$

$$\phi_m(x) = \begin{cases} 1 & \text{if } \bar{x} < \dfrac{\theta_0 + \theta_1}{2} \text{ if } \mu_1 < \mu_0 \\ 0 & \text{otherwise} \end{cases}$$

**Consider the case when $\sigma^2$ is unknown.**

For this case we propose a test which rejects $H_0$ when

$$\frac{\hat{L}(x/\theta_1)}{\hat{L}(x/\theta_0)} \geq 1$$

where $\hat{L}(x/\theta_i)$, (i=0,1) is the maximum of the likelihood under $H_i$ obtained from $L(x/\theta_i)$ by replacing $\sigma^2$ by its maximum likelihood estimate

$\widehat{\sigma_i^2} = \frac{1}{n}\sum_{j=1}^{n}\left(x_j - \theta_i\right)^2$; i=0,1.

Let $\phi_p(x)$ denote the critical function of the proposed test, then

$$\phi_p(x) = \begin{cases} 1 & if\ \hat{L}(x/\theta_1)\ \geq \hat{L}(x/\theta_0) \\ 0 & otherwise \end{cases}$$

On simplification we get

$$\phi_p(x) = \begin{cases} 1 & if\ \bar{x} \geq \frac{\theta_0+\theta_1}{2} \quad if\ \theta_1 > \theta_0 \\ o & otherwsie \end{cases}$$

and $$\phi_p(x) = \begin{cases} 1 & if\ \bar{x} < \frac{\theta_0+\theta_1}{2} \quad if\ \theta_1 < \theta_0 \\ o & otherwsie \end{cases}$$

Thus the proposed test $\phi_p(x)$ is equivalent to $\phi_m(x)$ which is the optimum test that minimizes the sum of two error probabilities $(\alpha + \beta)$. Thus we see that one gets the same test which minimises the sum of the two error probabilities irrespective of whether $\sigma^2$ is known or unknown.

**Acknowledgement** : Authors are thankful to Prof. Jokhan Singh for his guidance in writing this article.

**References**

Pandey, M. and Singh,J. (1978) : A note on testing mean of a normal population when variance is unknown. J.I.S.A. 16, pp. 145-146.

Pradhan,M. (1968) : On testing a simple hypothesis against a simple alternative making the sum of the probabilities of two errors minimum. J.I.S.A. 6, pp. 149-159.

Rao, C.R.(1973) : Linear Statistical Inference and its Applications. John Wiley and Sons

Singh,J. (1985): Optimum Statistical Test Procedures. Vikram Mathematical Journal Vol V, pp. 53-56.